\documentclass{amsart}
\usepackage{amsmath,amsthm,amsfonts,amssymb,hyperref}

\title{Existence of quasi-arcs}
\author[J. M. Mackay]{John M. Mackay}
\address{Department of Mathematics\\ University of Michigan\\
 Ann Arbor\\ Michigan 48109-1109}
\curraddr{Department of Mathematics\\ Yale University\\ New Haven\\
  Connecticut 06520-8283}
\email{jmmackay@umich.edu}
\thanks{This research was partially supported by NSF grant DMS-0701515.}
\date{December 10, 2007}

\newtheorem{thm}{Theorem}[section]
\newtheorem{cor}[thm]{Corollary}
\newtheorem{lem}[thm]{Lemma}
\newtheorem{prop}[thm]{Proposition}
\newtheorem*{claimdelta}{Claim $\Delta(n)$}

\theoremstyle{definition}

\theoremstyle{remark}

\DeclareMathOperator{\diam}{diam}
\newcommand{\eps}{\epsilon}
\newcommand{\del}{\delta}

\begin{document}

\begin{abstract}
We show that doubling, linearly connected metric spaces are quasi-arc
 connected.  This gives a new and short proof of a theorem of Tukia.
\end{abstract}

\maketitle

\section{Introduction}
\label{sec-intro}

It is a standard topological fact that a complete metric space
which is locally connected, connected and locally compact is
arc-wise connected.  Tukia~\cite{MR1420592} showed that an
analogous geometric statement is true: if a complete metric space
is linearly connected and doubling, then it is connected by
quasi-arcs, quantitatively.
In fact, he proved a stronger
 result: any arc in such a space may be approximated by a local quasi-arc
 in a uniform way.  In this note we give a new and more
direct proof of this fact.

This result is of interest in studying the quasisymmetric geometry of
metric spaces.  Such geometry arises in the study of the
boundaries of hyperbolic groups; Tukia's result was used in this
context by Bonk and Kleiner~\cite{MR2146190}, and also by the
author~\cite{Cdimpreprint}.
(Bonk and Kleiner use Assouad's embedding theorem to translate Tukia's
result from its original context of subsets of $\mathbb{R}^n$
into our setting of doubling and linearly connected metric spaces.)

Before stating the theorem precisely, we recall some definitions.
A metric space $(X,d)$ is said to be {\em doubling} if
there exists a constant $N$ such that every ball
can be covered by at most $N$ balls of half the radius.
Note that any complete, doubling metric space is proper:
all closed balls are compact.

We say $(X,d)$  is {\em $L$-linearly connected} for some $L \geq 1$ if
 for all $x,\,y \in X$ there exists a compact, connected set $J \ni x,\,y$
 of diameter less than or equal to $L d(x,y)$.
(This is also known as bounded turning or LLC(1).)
We can actually assume that $J$ is an arc, at the
cost of increasing $L$ by an arbitrarily small amount.
To see this, note that $X$ is locally connected, and
so the connected components of an open set are open.
Thus, for any open neighborhood $U$ of $J$,
the connected component of $U$ that contains $J$ is an open set.
 We can replace $J$ inside $U$ by an arc with the same endpoints,
since any open, connected subset of a locally compact, locally connected
metric space is arc-wise connected~\cite[Corollary 32.36]{MR0221455}.

For any $x$ and $y$ in an embedded arc $A$, we denote by $A[x,y]$
the closed, possibly trivial, subarc of $A$ that lies between them.
We say that an arc $A$ in a doubling and complete metric space is an
{\em $\eps$-local $\lambda$-quasi-arc} if
$\diam(A[x,y]) \leq \lambda d(x,y)$ for all $x,\, y \in A$ such that
$d(x,y) \leq \eps$.
(This terminology is explained by Tukia and V\"ais\"al\"a's
 characterization of quasisymmetric images of the unit interval
 as those metric arcs that are doubling and bounded turning~\cite{MR595180}.)

One non-standard definition will be useful in our exposition.
We say that an arc $B$ {\em $\eps$-follows} an arc $A$ if there
exists a coarse map $p:B \rightarrow A$, sending endpoints to endpoints,
such that for all $x,\,y \in B$, $B[x,y]$ is in the $\eps$-neighborhood of
$A[p(x),p(y)]$; in particular, $p$ displaces points at most $\eps$.
(We call the map $p$ {\em coarse} to emphasize that it is not
 necessarily continuous.)

The condition that $B$ $\eps$-follows $A$ is stronger than the
condition that $B$ is contained in the $\eps$-neighborhood of $A$.
It says that, coarsely, the arc $B$ can be obtained from the arc $A$
by cutting out `loops.'  (Of course, $A$ contains no
 actual loops, but it may have subarcs of large diameter whose endpoints
 are $2\eps$-close.)

We can now  state the stronger version of Tukia's theorem precisely, and
as an immediate corollary our initial statement
\cite[Theorem 1B, Theorem 1A]{MR1420592}:

\begin{thm}[Tukia]\label{thm-tukia}
 Suppose $(X,d)$ is a L-linearly connected,  N-doubling,
 complete metric space.
 For every arc $A$ in $X$ and every $\eps >0$, there is an arc $J$
 that $\eps$-follows $A$, has the same endpoints as $A$,
 and is an $\alpha\eps$-local $\lambda$-quasi-arc,
 where $\lambda = \lambda(L,N) \geq 1$ and $\alpha = \alpha(L,N) >0$.
\end{thm}

\begin{cor}[Tukia]\label{cor-tukia}
 Every pair of points in a $L$-linearly connected, $N$-doubling,
  complete metric space is
 connected by a $\lambda$-quasi-arc, where $\lambda = \lambda(L,N) \geq 1$.
\end{cor}

Our strategy for proving Theorem~\ref{thm-tukia} is straightforward:
 find a method of straightening an arc on a given
 scale (Proposition~\ref{lem-coarsestr}), then
 apply this result on a geometrically decreasing sequence
 of scales to get the desired local quasi-arc as a limiting object.
 The statement of this proposition and the resulting proof of the theorem
 essentially follow Tukia~\cite{MR1420592}, but the
 proof of the proposition is new and much shorter.  We include a
 complete proof for convenience to the reader.

The author thanks Mario Bonk and, in particular,
his advisor Bruce Kleiner for many
helpful suggestions and fruitful conversations.

\section{Main Results}
\label{sec-tukiaproof}

Given any arc $A$ and $\iota >0$, the following
 proposition allows us to straight\-en $A$ on a scale $\iota$ inside the
 $\iota$-neighborhood of $A$.

\begin{prop}\label{lem-coarsestr}
 Given a complete metric space $X$ that is
 $L$-linearly connected and $N$-doubling, there exist constants
 $s=s(L,N)>0$ and $S=S(L,N)>0$ with the following property:
 for each $\iota > 0$ and each arc $A \subset X$, there exists an arc
 $J$ that $\iota$-follows $A$, has the same endpoints as $A$,
 and satisfies
 \begin{equation} \label{eq-cqa}
  \forall x,y \in J,\  d(x,y) < s\iota \implies
   \diam(J[x,y]) < S\iota \tag{$\ast$}.
 \end{equation}
\end{prop}

We will apply this proposition on a decreasing sequence
of scales to get a local quasi-arc in the limit.  The
key step in proving this is given by the following
lemma.

\begin{lem} \label{lem-approx}
 Suppose $(X,d)$ is a L-linearly connected, N-doubling, complete
 metric space, and let  $s,\, S,\, \eps$ and $\del$
 be fixed positive constants satisfying
 $\del \leq \min\{\frac{s}{4+2S},\frac{1}{10}\}$.
 Now, if we have a sequence of arcs $J_1, J_2, \ldots, J_n, \ldots$ in $X$,
  such that for every $n \geq 1$
 \begin{itemize}
  \item $J_{n+1}$ $\eps \del^n$-follows $J_n$, and
  \item $J_{n+1}$ satisfies $($\ref{eq-cqa}$)$ with
   $\iota = \eps \del^n$ and $s,\,S$ as fixed above,
 \end{itemize}
 then the Hausdorff limit $J = \lim_\mathcal{H} J_n$
 exists, and is an $\eps \del^2$-local
 $\frac{4S+3\del}{\del^2}$-quasi-arc.

 Moreover, the endpoints of $J_n$ converge to the endpoints of $J$,
 and $J$ $\eps$-follows $J_1$.
\end{lem}

We shall need some standard definitions.  The (infimal) distance between two subsets
$U,\,V \subset X$ is defined as $d(U,V) = \inf\{d(u,v) : u \in U, v \in V\}$.
If $U=\{u\}$, then we set $d(u,V) = d(U,V)$.

The $r$-neighborhood of $U$ is the set $N(U,r) = \{ x : d(x,U) < r \}$,
and the Hausdorff distance between $U$ and $V$, $d_\mathcal{H}(U,V)$,
 is defined to be the infimal $r$ such that $U \subset N(V,r)$ and $V \subset N(U,r)$.
For more information, see~\cite[Chapter 7]{MR1835418}.

We will now prove Theorem~\ref{thm-tukia}.

\begin{proof}[Proof of Theorem~\ref{thm-tukia}]
 Let $s$ and $S$ be given by Proposition~\ref{lem-coarsestr},
 and set $\del = \min\{\frac{s}{4+2S},\frac{1}{10}\}$.

 Let $J_1 = A$ and apply Proposition~\ref{lem-coarsestr}
 to $J_1$ and $\iota = \eps \del$ to get an arc $J_2$
 that $\eps \del$-follows $J_1$.
   Repeat, applying the lemma to $J_n$ and
 $\iota = \eps \del^n$,
 to get a sequence of arcs $J_n$, where each $J_{n+1}$
 $\eps \del^n$-follows $J_n$,
 and satisfies (\ref{eq-cqa}) with $\iota = \eps \del^n$.

 We can now apply Lemma~\ref{lem-approx} to find an
 $\alpha\eps$-local $\lambda$-quasi-arc $J$ that
 $\eps$-follows $A$,
 where $\alpha = \del^2$ and $\lambda = \frac{4S+3\del}{\del^2}$.
 Every $J_n$ has the same endpoints as $A$,
 so $J$ will also have the same endpoints.
\end{proof}

The proof of Lemma~\ref{lem-approx} relies on some fairly
straightforward estimates and a classical characterization
of an arc.

\begin{proof}[Proof of Lemma~\ref{lem-approx}]
 For every $n \geq 1$, $J_{n+1}$ $\eps\del^n$-follows $J_n$.
 We denote the associated coarse map by $p_{n+1}:J_{n+1} \to J_n$.

 In the following, we will make frequent use
 of the inequality
 $\sum_{n=0}^{\infty} \del^n < \frac{11}{9}$.

 We begin by showing that the Hausdorff limit
 $J = \lim_{\mathcal{H}} J_n$ exists.
 The collection of all compact subsets of a compact metric space,
 given the Hausdorff metric, is itself a compact metric
 space~\cite[Theorem 7.3.8]{MR1835418}.
 Since $\{ J_n \}$ is a sequence of compact sets in a
 bounded region of a proper metric space, to show that
 the sequence converges with respect to the Hausdorff metric,
 it suffices to show that the sequence is Cauchy.

 One bound follows by construction:
  $J_{n+m} \subset N(J_n, \frac{11}{9}\eps\del^n)$ for all $m \geq 0$.
 For the second bound, take $J_{n+m}$ and split it
 into subarcs of diameter at most $\eps\del^n$, and write this as
 $J_{n+m} = J_{n+m}[z_0,z_1] \cup \cdots \cup J_{n+m}[z_{k-1},z_k]$ for some
 $z_0, \ldots, z_k$ and some $k >0$.
 Our coarse maps compose to give $p:J_{n+m} \rightarrow J_n$,
 showing that $J_{n+m}$ $\frac{11}{9} \eps\del^n$-follows $J_n$.
 Furthermore, since
 $d(z_i, z_{i+1}) \leq \eps\del^n$, we have
  $d(p(z_i), p(z_{i+1})) \leq 4\eps\del^n \leq s\eps\del^{n-1}$.
 Combining this with the fact that $p$ maps endpoints to endpoints, for $n \geq 2$
 we have
\begin{align*}
  J_n & = J_n[p(z_0),p(z_1)] \cup \cdots \cup J_n[p(z_{k-1}),p(z_k)]
	\subset N(\{p(z_0), \ldots, p(z_k)\}, S\eps\del^{n-1}) \\
      & \subset N\left(J_{n+m}, \frac{11}{9}\eps\del^n + S\eps\del^{n-1}\right).
\end{align*}

 Taken together, these bounds give
 $d_\mathcal{H}(J_n,J_{n+m}) \leq \frac{11}{9}\eps\del^n + S\eps\del^{n-1}$,
   so $\{ J_n \}$ is Cauchy and the limit $J= \lim_{\mathcal{H}} J_n$ exists.
 Moreover, $J$ is compact (by definition) and connected (because
 each $J_n$ is connected).

 Now we let $a_n$, $b_n$ denote the endpoints of $J_n$.
 Since $p_n(a_n) = a_{n-1}$, and $p_n$ displaces
 points at most $\eps \del^n$, the sequence $\{a_n\}$
 is Cauchy and hence converges to some point $a \in J$.
 Similarly, $\{b_n\}$ converges to a point $b \in J$.

 There are two cases to consider.
 If $a = b$, then $d(a_n, b_n) \leq 2 \frac{11}{9}\eps\del^n
  \leq s\eps\del^{n-1}$.
 Consequently, $\diam(J_n) \leq S\eps\del^{n-1}$,
 $J = \lim_\mathcal{H} J_n$ has diameter zero, and thus $J= \{a\}$.
 Otherwise, $a \neq b$ and so $J$ is non-trivial.  We claim that
 in this case $J$ is a local quasi-arc.

 To show $J$ is an arc with endpoints $a$ and $b$
 it suffices to demonstrate that every point $x \in J \setminus \{a,b\}$
 is a cut point \cite[Theorems 2-18 and 2-27]{MR0125557}.
 The topology of $J_n$ induces an order on $J_n$
 with least element $a_n$ and greatest $b_n$.
 Given $x \in J$, we define three points $h_n(x)$, $x_n$ and $t_n(x)$ that
 satisfy $a_n < h_n(x) < x_n < t_n(x) < b_n$, where $x_n$ is chosen such that
 $d(x,x_n) \leq \frac{11}{9}\eps\del^n$, and $h_n(x)$ and $t_n(x)$ are the
 first and last elements of $J_n$ at distance $(S+1)\eps\del^{n-1}$ from $x$.
 As long as $x$ is not equal to $a$ or $b$, for $n$ greater than some $n_0$
 these points will exist and this definition will be valid.

 We shall denote the $\frac{11}{9}\eps\del^n$-neighborhoods
 of $J_n[a_n,h_n(x)]$ and $J_n[t_n(x),b_n]$ by $H_n(x)$ and $T_n(x)$ respectively,
 and define
 $H(x) = \cup \{H_n(x) : n \geq n_0\}$ (the Head)
 and $T(x) = \cup \{T_n(x) : n \geq n_0\}$ (the Tail).
 By definition, $H(x)$ and $T(x)$ are open.
 We claim that, in addition, they are disjoint and
 cover $J \setminus \{x\}$, and so $x$ is a cut point.

 Fix $y \in J$, and suppose $y \notin H(x) \cup T(x)$.
 We want to show that $y = x$.
 To this end, we bound the diameter of $J_n[h_n(x),t_n(x)]$ using $J_{n-1}$.
 Because $d(p_{n}(h_n(x)),p_{n}(t_n(x))) \leq
   2\eps\del^{n-1} +
   2(S+1)\eps\del^{n-1} \leq s\eps\del^{n-2}$,
 we know that the diameter of
 $J_{n-1}[p_{n}(h_n(x)),p_{n}(t_n(x))]$ must be less than $S\eps\del^{n-2}$.
 Thus the diameter of $J_n[h_n(x),t_n(x)]$ is less than
  $S\eps\del^{n-2}+2\eps\del^{n-1}$, as $J_n$
  $\eps\del^{n-1}$-follows $J_{n-1}$.

 For every $n \geq n_0$, $y$ is $\frac{11}{9}\eps\del^n$ close to some $y_n \in J_n$.
 Since $y \notin H(x) \cup T(x)$, $y_n$ must lie in $J_n[h_n(x),t_n(x)]$, so
 \begin{align*}
  d(x,y) & \leq d(x,J_n[h_n(x),t_n(x)]) + \diam(J_n[h_n(x),t_n(x)])
     + d(y_n,y) \\
   & \leq 2\frac{11}{9}\eps\del^n + (S+2\del)\eps\del^{n-2}
    = \left(2\frac{11}{9}\del^2+S+2\del \right)\eps\del^{n-2},
 \end{align*}
 therefore $d(x,y) = 0$ and $J \setminus (H(x) \cup T(x)) = \{x\}$.

 We now show that $H(x)$ and $T(x)$ are disjoint.
 If not, then $H_n(x) \cap T_m(x) \neq \emptyset$ for some $n$ and $m$.
 It suffices to assume $n \leq m$.  Now
  $T_m(x) \subset N(J_m[x_m,b_m],\frac{11}{9}\eps\del^m)$ by definition.
 We send $J_m$ to $J_n$ using
 $f = p_{n+1} \circ \cdots \circ p_m : J_m \rightarrow J_n$, to get that
 $T_m(x) \subset N(J_n[f(x_m),b_n],3\eps\del^n)$.
 Since
 \[ d(f(x_m),x_n) \leq d(f(x_m),x_m)+d(x_m,x) + d(x,x_n)
   < 4\eps\del^n < s\eps\del^{n-1} \]
 we have, even for $n=m$,
 \[ T_m(x) \subset N(J_n[x_n,b_n],3\eps\del^n) \cup
   B(x_n, (S+3\del)\eps\del^{n-1}) .\]

 Since $(S+3\del)\eps\del^{n-1} + \frac{11}{9}\eps\del^n
  < (S+\frac{1}{2})\eps\del^{n-1}$, $H_n(x)$ cannot meet $T_m(x)$
  in the ball $B(x_n, (S+3\del)\eps\del^{n-1})$.
 Thus $H_n(x) \cap T_m(x) \neq \emptyset$ implies that there exist
 points $p$ and $q$ in $J_n$ such that
   $a_n \leq p \leq h_n(x) < x_n \leq q \leq b_n$ and
   $d(p,q) < 3\eps\del^n < s\eps\del^{n-1}$.
 But then we know that $J_n[p,q]$ has diameter less than
 $S\eps\del^{n-1}$, while containing both $h_n(x)$ and $x_n$.
 This contradicts the definition of $h_n(x)$, so $H(x) \cap T(x) = \emptyset$.

 We have shown that $J$ is an arc with endpoints $a$ and $b$;
 it remains to show that $J$ is a local quasi-arc with the
 required constants.

 Say we are given
 $x$ and $y$ in $J$, with $x_n$ and $y_n$ as before.
 Our argument shows that the segments $J_n[x_n,y_n]$ converge to some arc
 $\tilde{J}[x,y]$, because $J_{n+1}[x_{n+1},y_{n+1}]$
  $(\eps \del^n +S\eps\del^{n-1})$-follows $J_n[x_n,y_n]$ for
 all $n \geq 2$.
 This arc $\tilde{J}[x,y]$ must lie in $J$,
 therefore $\tilde{J}[x,y]$ must equal $J[x,y]$.
 Now, suppose that
   $d(x,y) \in (\eps\del^{n+1},\eps\del^n]$
   holds for some $n\geq 2$.
 Then
 $d(x_n,y_n) \leq 3\eps\del^n + \eps\del^n < s\eps\del^{n-1}$,
 and so the subarc $J[x,y]$, which lies in
  $N(J_n[x_n,y_n],\frac{11}{9}\eps(\del^n+S\del^{n-1}))$,
 has diameter less than
  $S\eps\del^{n-1}+3\eps(\del^n+S\del^{n-1}) \leq
    \frac{4S+3\del}{\del^2}d(x,y)$, as desired.

 Furthermore, this same argument gives that, for all $n \geq 2$, $J$
 $\frac{11}{9}\eps(\del^n+S\del^{n-1})$-follows $J_n$,
 which itself $\frac{11}{9}\eps\del$-follows $J_1=A$.
 Taking $n$ sufficiently large, we have that
 $J$ $\eps$-follows $A$.
\end{proof}

\section{Discrete paths and the proof of Proposition~\ref{lem-coarsestr}}
\label{sec-mainprop}

The proof of Proposition~\ref{lem-coarsestr} is based on a quantitative
 version of a simple geometric result.
Before we state this result, recall that a maximal $r$-separated set $\mathcal{N}$ is a subset of $X$
  such that for all distinct $x,\,y \in \mathcal{N}$ we have $d(x,y) \geq r$, and
  for all $z \in X$ there exists some $x \in \mathcal{N}$ with $d(z,x) < r$.

Now suppose that we are given a maximal
 $r$-separated set $\mathcal{N}$ in an $L$-linearly connected,
 $N$-doubling, complete metric space $X$.  Then we can find
 a collection of sets $\{V_x\}_{x \in \mathcal{N}}$
 so that each $V_x$ is a connected union of finitely many arcs in $X$,
 and for all $x,\, y \in \mathcal{N}$:

 \begin{enumerate}
  \item[(1)] $d(x,y) \leq 2r \implies y \in V_x$.
  \item[(2)] $\diam(V_x) \leq 5Lr$.
  \item[(3)] $V_x \cap V_y = \emptyset \implies d(V_x, V_y)  > 0$.
 \end{enumerate}

 For $x \in \mathcal{N}$, we can construct each $V_x$ by defining it to be the union of
  finitely many arcs joining $x$ to each $y \in \mathcal{N}$ with
  $d(x,y) \leq 2r$.  By linear connectedness, we can ensure that $\diam(V_x) \leq 4Lr$.
  Condition (3) is trivially satisfied for compact subsets of a metric
  space, but we will strengthen it to the following:

 \begin{enumerate}
  \item[(3$^\prime$)] $V_x \cap V_y = \emptyset \implies d(V_x, V_y)  > \del r$.
 \end{enumerate}

\begin{lem}\label{lem-nicenet}
 We can construct the sets $V_x$ satisfying (1), (2) and (3$\,^\prime\!$) for
 $\del = \del(L,N)$.
\end{lem}
\begin{proof}
  Without loss of generality, we can rescale the metric to set $r=1$.

  Since $X$ is doubling, the maximum number of 1-separated points
  in a $20L$-ball is bounded by a constant $M=M(20L,N)$.
  We can label every point of $\mathcal{N}$ with
  an integer between $1$ and $M$, such that no two points
  have the same label if they are separated by a distance less than $20L$.

  To find this labelling, consider the collection of all
  such labellings on subsets of $\mathcal{N}$ under
  the natural partial order. A Zorn's Lemma argument gives
  the existence of a maximal element: our desired labelling.
  So $\mathcal{N}$ is the disjoint union of $20L$-separated sets
  $\mathcal{N}_1, \ldots, \mathcal{N}_M$.

  Now let $\mathcal{N}_{\leq n} = \cup_{k=1}^{n} \mathcal{N}_k$,
  and consider the following

\begin{claimdelta}
  We can find $V_x$
  for all $x \in \mathcal{N}_{\leq n}$, such that
  for all $x, y \in \mathcal{N}_{\leq n}$ (1), (2) and (3$\,^\prime\!$)
  are satisfied with $\del = \frac{1}{2}(5L)^{-n}$.
\end{claimdelta}

  $\Delta(0)$ holds trivially, and Lemma~\ref{lem-nicenet} immediately
  follows from $\Delta(M)$, with $\del = \del(L,N) = \frac{1}{2}(5L)^{-M}$.
  So we are done, modulo the statement that $\Delta(n) \implies \Delta(n+1)$
  for $n < M$.
\end{proof}

\begin{proof}[Proof that $\Delta(n) \implies \Delta(n+1)$, for $n < M$]
 By $\Delta(n)$, we have sets $V_x$ for all $x$ in $\mathcal{N}_{\leq n}$.

  As $\mathcal{N}_{n+1}$ is $20L$-separated we can treat the constructions
  of $V_x$ for each $x \in \mathcal{N}_{n+1}$ independently.
  We begin by creating a set $V_x^{(0)}$ that is the union of
  finitely many arcs joining $x$ to its
  $2$-neighbors in $\mathcal{N}$.  We can ensure that $\diam(V_x^{(0)}) \leq 4L$.

  Now construct $V_x^{(i)}$ inductively, for $1 \leq i \leq n$.
  $V_x^{(i-1)}$ can be $5L$-close to at most one $y \in \mathcal{N}_i$.
  If $d(V_x^{(i-1)}, V_y) \in (0, \frac{1}{2} (5L)^{-i})$, then define
  $V_x^{(i)}$ by adding to $V_x^{(i-1)}$ an arc of diameter at most
  $ L (5L)^{-i}$ joining $V_x^{(i)}$ to $V_y$.
  Otherwise, let $V_x^{(i)} = V_x^{(i-1)}$.  Continue until $i=n$
  and set $V_x = V_x^{(n)}$.

  Note that $V_x$ satisfies (1) and (2) by construction.  The
  only non-trivial case we need to check for (3$^{\prime}$) is whether
  $d(V_x, V_y) \in (0, \frac{1}{2} (5L)^{-n})$ for some
  $y \in \mathcal{N}_i$, $i < n$.  (The $i=n$ case follows from the
  last step of the construction.)
  Then, since $V_x = V_x^{(n)} \supset V_x^{(i)}$,
  $V_x^{(i)} \cap V_y \neq \emptyset$, and
  $d(V_x^{(i)}, V_y) \geq \frac{1}{2} (5L)^{-i}$.
  The construction then implies that
  \begin{align*}
     d(V_x, V_y) & \geq \frac{1}{2} (5L)^{-i}( 1 - (2L)(5L)^{-1} -
           (2L)(5L)^{-2} - \cdots - (2L)(5L)^{-(n-i)}) \\
         & > \frac{1}{2} (5L)^{-n} (5L) \left(1 - \frac{2/5}{1-(1/(5L))} \right)
        \geq \frac{5}{2} \left( \frac{1}{2} (5L)^{-n} \right),
  \end{align*}
  contradicting our assumption, so $\Delta(n+1)$ holds.
\end{proof}

  We now finish by using this construction to prove our proposition.

\begin{proof}[Proof of Proposition~\ref{lem-coarsestr}]
  By rescaling the metric, we may assume that $\iota = 20L$.
  If $d(a,b) \leq 20 = \frac{\iota}{L}$, then join $a$ to
  $b$ by an arc of diameter less than $\iota$.  This arc
  will, trivially, satisfy our conclusion for any $S \geq 1$.

  Otherwise, $d(a,b) > 20$.  In the hypotheses for Lemma~\ref{lem-nicenet},
  let $r=1$ and let $\mathcal{N}$
  be a maximal $1$-separated set in $X$ that contains both $a$ and $b$.
  Now apply Lemma~\ref{lem-nicenet} to get $\{V_x\}_{x \in \mathcal{N}}$
  satisfying (1), (2) and (3$^\prime$) for $\del = \del(L,N) > 0$.

  We want to `discretize' $A$ by finding a corresponding sequence of
  points in $\mathcal{N}$.  Now, every open ball in $X$ meets the
  arc $A$ in a collection of disjoint, relatively open intervals.
  Since $\mathcal{N}$ is a maximal
  $1$-separated set, the collection of open balls
  $\{ B(x, 1) : x \in \mathcal{N}\}$ covers $X$; in particular, it
  covers $A$.  By the compactness of $A$, we can find a finite cover
  of $A$ by connected, relatively open intervals, each lying
  in some $B(x,1)$, $x \in \mathcal{N}$.

  Using this finite cover, we can find
  points $x_i \in \mathcal{N}$ and $y_i \in A$ for $0 \leq i \leq n$,
  such that $a=y_0 < \cdots < y_n=b$ in the order on $A$,
  and $A[y_i,y_{i+1}] \subset B(x_i, 1)$ for each $0 \leq i < n$.
  Since $a,\,b \in \mathcal{N}$, we have that $x_0 = a$ and $x_n = b$.
  The sequence $(x_0, \ldots, x_n)$ is a discrete path in $\mathcal{N}$
  that corresponds naturally to $A$.

  We now find a subsequence $(x_{r_j})$ of $(x_i)$ such that the
  corresponding sequence of sets $(V_{x_{r_j}})$ forms a `path'
  without repeats.
  Let $r_0 = 0$, and
  for $j \in \mathbb{N}^+$ define $r_j$ inductively as
    $r_j = \max\{ k : V_{x_k} \cap V_{x_{r_{j-1}}} \neq \emptyset\}$,
    until $r_m = n$ for some $m \leq n$.
  The integer $r_j$ is well defined
  since $d(y_{(r_{j-1}+1)},x_k) \leq 1$ for $k=r_{j-1}$ and $k=r_{j-1}+1$,
  so $V_{x_{(r_{j-1}+1)}} \cap V_{x_{r_{j-1}}} \neq \emptyset$.
  Note that if $i + 2 \leq j$ then $V_{x_{r_i}} \cap V_{x_{r_j}} = \emptyset$,
  and thus $d(V_{x_{r_i}}, V_{x_{r_j}}) > \del$.

  Let us construct our arc $J$ in segments.  First, let $z_0 = x_{r_0}$.
  Second, for each $i$ from $0$ to $m-1$,
  let $J_i = J_i[z_i, z_{i+1}]$ be an arc in $V_{x_{r_i}}$ that
  joins $z_i \in V_{x_{r_i}}$
  to some $z_{i+1} \in V_{x_{r_{i+1}}}$, where
  $z_{i+1}$ is the first point of $J_i$ to meet
  $V_{x_{r_{i+1}}}$.  (In the case $i=m-1$, join $z_{m-1}$ to
  $x_{r_m} = z_m$.)  Set $J = J_0 \cup \cdots \cup J_m$.

  This path $J$ is an arc since each $J_i \subset V_{x_{r_i}}$ is an arc, and
  if there exists a point $p \in J_i \cap J_j$ for some $i<j$, then
  $j = i+1$ and $p = z_{i+1} = z_j$.
  This is true because
  $V_{x_{r_i}} \cap V_{x_{r_j}} \neq \emptyset$ implies that $j = i+1$,
  and the definition of $z_{i+1}$ implies that
  $J_i \cap V_{x_{r_{i+1}}} = \{z_{i+1}\}$.
  Any finite sequence of arcs that meet only at consecutive endpoints is also
  an arc, so we have that $J$ is an arc.

  In fact, $J$ satisfies (\ref{eq-cqa}).
  For any $y$, $y' \in J$, $y<y'$, we can
  find $i \leq j$ such that $z_i \leq y < z_{i+1}$, $z_j \leq y' < z_{j+1}$.
  (If $y = z_m$, set $i=m$; likewise for $y'$.)
  If $d(y,y') \leq \del$ then,
  because $y \in V_{x_{r_i}}$ and $y' \in V_{x_{r_j}}$, we have
  $d(V_{x_{r_i}}, V_{x_{r_j}}) \leq \del$, so either $j=i$ or $j=i+1$.
  This gives that $J[y,y'] \subset V_{x_{r_i}} \cup V_{x_{r_j}}$, and so
  $\diam(J[y,y'])$ is bounded above by $10L$.

  Furthermore, $J$ $\iota$-follows $A$.
  There is a coarse map $f:J \rightarrow A$ defined by
  the following composition: first map $J$ to $\mathcal{N}$ by
   sending $y \in J[z_i,z_{i+1}) \subset J$
   to $x_{r_i} \in \mathcal{N}$, and sending $x_{r_m}$ to itself.
  Second, map each $x_{r_i}$ to the
  corresponding $y_{r_i}$ in $A$.
  Taking arbitrary $y < y'$ in $J$ as before, we see that
  \begin{align*}
   J[y,y'] & \subset J[z_i,z_{j+1}] \subset N(\{x_{r_i}, \ldots, x_{r_j}\},5L)
    \subset N(\{y_{r_i}, \ldots, y_{r_j}\},5L+1) \\
    & \subset N(A[y_{r_i},y_{r_j}],5L+1)
    \subset N(A[f(y),f(y')],\iota).
  \end{align*}

  We let $s = \frac{1}{20L}\del$ and $S = \frac{1}{20L}10L$, and have
  proven the Proposition.
\end{proof}

{\em Remark:}  This method of proof allows one to explicitly estimate the
constants given in the statements of Theorem~\ref{thm-tukia} and
Corollary~\ref{cor-tukia}, but for most applications this is not
necessary.


\bibliographystyle{amsplain}
\bibliography{confdim}

\end{document}